%% file: GAforEGFinal.tex
\documentclass{birkjour}
\def\isFV{false}
\def\isBeamer{false}
\input{GunnMac}

\usepackage[title]{appendix}
\let\OLDthebibliography\thebibliography
\renewcommand\thebibliography[1]{
  \OLDthebibliography{#1}
  \setlength{\parskip}{0pt}
  \setlength{\itemsep}{0pt plus 0.3ex}
}

\graphicspath{{./pictures/}{./}{./defense/pictures/}}

\renewcommand{\vec}[1]{\mathbf{#1}}

\newcommand{\gTh}{\cite{gunnThesis}\xspace}

\newcommand{\thetask}{{the chosen task}\xspace}
\begin{document}
\title{Geometric algebras for euclidean geometry\footnote{This article has been published as \cite{gunn2016}.  The final publication is available at link.springer.com.}}
\ifthenelse{\equal{\isFV}{true}}{
}
{
\author{Charles Gunn}
\address{Instit\"{u}t fŸr Mathematik  MA 8-3\\
Technische Universit\"{a}t Berlin\\
Str. des 17 Juni 136\\
10623 Berlin Germany\\
\\Charles Gunn\\
Raum+Gegenraum \\  Brieselanger Weg 1 \\ 14612 Falkensee  }
\email{gunn@math.tu-berlin.de}
}
\maketitle

\begin{abstract}
The discussion of how to apply geometric algebra to euclidean $n$-space has been clouded by a number of conceptual misunderstandings which we first identify and resolve, based on a thorough review of crucial but largely forgotten themes from $19^{th}$ century mathematics.
We then introduce the dual projectivized Clifford algebra $\pdclal{n}{0}{1}$ (euclidean PGA) as the most promising homogeneous (1-up) candidate for euclidean geometry.   We compare  euclidean PGA and the popular 2-up model CGA (conformal geometric algebra),  restricting attention to flat geometric primitives, and show that on this domain they exhibit the same formal feature set. We thereby establish that euclidean PGA is the smallest structure-preserving euclidean GA.  We compare the two algebras in more detail, with respect to a number of practical criteria, including implementation of kinematics and rigid body mechanics.  We then extend the comparison to include euclidean sphere primitives.  We conclude that euclidean PGA provides a natural transition,  both scientifically and pedagogically, between vector space models and the more complex and powerful CGA.  
\end{abstract}

\section{Introduction}

Although noneuclidean geometry of various sorts plays a fundamental role in theoretical physics and cosmology, the overwhelming volume of practical science and engineering takes place  within classical euclidean space $\Euc{n}$.  For this reason it is of no small interest to establish the best computational model for this space.    
In particular, this article explores the question, which form of geometric algebra is best-suited for computing in euclidean space?  In order to have a well-defined domain for comparison,  we restrict ourselves for most of the article to flat geometric primitives (points, lines, planes, and their higher-dimensional analogs), along with kinematics and rigid body mechanics.   We refer to this domain  as \quot{the chosen task}. In \Sec{sec:round} we widen the primitive set to include spheres.

\subsection{Overview}
As the later results of the article are based upon  crucial but often overlooked $19^{th}$ century mathematics, the article begins with the latter.  In \Sec{sec:qbrbm} we review the quaternions and the biquaternions, focusing on their use to model euclidean kinematics and rigid body mechanics, a central part of \thetask.  
In \Sec{sec:featga} we trace the path from  exterior algebra to geometric algebra in the context of projective geometry, paying special attention to the dual exterior algebra and the Cayley-Klein construction of metric spaces within projective space.   This culminates with the introduction of projective geometric algebra (PGA), a homogeneous model for euclidean (and other constant-curvature metric) geometry.   In \Sec{sec:clar} we demonstrate that PGA is superior to other proposed homogeneous models.  \Sec{sec:hmudm} discusses and disposes of a variety of  misconceptions appearing in the geometric algebra literature regarding geometric algebras with degenerate metrics (such as PGA). 
In \Sec{sec:comp} we identify a feature set for a euclidean geometric algebra for \thetask, and verify that both PGA and CGA fulfil all features.  \Sec{sec:comparison} turns to a comparison based on a variety of practical criteria, such as the implementation of kinematics and rigid body mechanics.  \Sec{sec:round} extends the comparison to include spheres, and concludes that the \quot{roundness} of CGA has both positive and negative aspects.  Finally, \Sec{sec:conclusion} positions PGA as a natural stepping stone, both scientifically and pedagogically, between vector space geometric algebra and CGA.  

\section{Quaternions, Biquaternions, and Rigid Body Mechanics}
\label{sec:qbrbm}

The quaternions (Hamilton, 1844) and the biquaternions (\cite{clifford73}) are important forerunners of  geometric algebra.  They exhibit most of the important features of a geometric algebra, such as an associative geometric product consisting of a symmetric (\quot{inner}) part and an anti-symmetric (\quot{outer}) part, and the ability to represent isometries as sandwich operators.  However,  they lack the graded algebra structure possessed by geometric algebra.   

They have a special importance in our context, since the even subalgebra of PGA (see below, \Sec{sec:pga}) is, for $n=3$ (the case of most practical interest), isomorphic to the biquaternions.  \emph{All the desirable features of the biquaternions are then inherited by PGA}.  In particular, the biquaternions contain a model of kinematics and rigid body mechanics, an important component of \thetask, that  compares favorably with modern alternatives (see \Sec{sec:krbdasct}). The biquaternions  also reappear below in \Sec{sec:lackcov}.
Because of this close connection with the themes of this article, and the absence of a comparable treatment in the literature, we give a brief formulation of relevant results here.  We assume the reader has an introductory acquaintance with quaternions and biquaternions.

\subsection{Quaternions} 
We first show how the Euler top can be advantageously represented using  quaternions.  
A more detailed treatment is available in Ch. 1 of \gTh.
Let $\mathbb{UH}$ represent the unit quaternions, and $\mathbb{IH}$ represent the purely imaginary quaternions.
Every unit quaternion can be written as the exponential of an imaginary quaternion.  A rotation $\vec{R}$ around the  unit vector $(b,c,d)$ through an angle $\alpha$ can be represented by a \quot{sandwich operator} as follows. Define the imaginary quaternion $\vec{m} := b\vec{i} + c\vec{j} + d\vec{k}$ which then gives rise to the unit quaternion $\vec{g} := e^{\frac{\alpha}{2}\vec{m}}$ by exponentiation.  Let $\vec{x} := x\vec{i} + y\vec{j} + z\vec{k}$ represent an arbitrary vector $(x,y,z) \in \R{3}$.  Then the rotation $\vec{R}$ applied to $(x,y,z)$ is given in the quaternion product as $\vec{x}' =  \vec{g} \vec{x} {\vec{g}^{-1}} (=  \vec{g} \vec{x} \overline{\vec{g}})$.  

\subsubsection{Quaternion model of Euler top} \label{sec:qmet} Recall that the Euler top is a rigid body constrained to move around its centre of gravity. We assume there are no  external forces. Let $\vec{g}(t)$, a path in $\mathbb{UH}$, be the motion of the rigid body.  Let $\vec{M}_{c} \in \mathbb{IH}$ and $\vec{V}_{c} \in \mathbb{IH}$ represent the instantaneous momentum and velocity, resp., in body coordinates. Considered as vectors in $\R{3}$, they are related by the inertia tensor $\vec{A}$ via $\vec{M}_{c} = \vec{A} ( \vec{V}_{c})$. Then the Euler equations of motion can be written using the quaternion product:
\begin{align*}
\dot{\vec{g}} &= \vec{g} \vec{V}_{c} \\
\dot{\vec{M}}_{c} &=  \frac12(\vec{V}_{c}\vec{M}_{c} - \vec{M}_{c} \vec{V}_{c})
\end{align*}
Notice that this representation has practical advantages over the traditional linear algebra approach using matrices: normalizing a quaternion brings it directly onto the 3D solution space of the unit quaternions, while the matrix group $SO(3)$ is a 3-dimensional subspace of the 9-dimensional space of 3x3 matrices.  Numerical integration proceeds much more efficiently in the former case since, after normalising, there are no chances for \quot{wandering off}, while the latter has a 6D space of invalid directions that lead away from $SO(3)$.  We meet the same problematic below in \Sec{sec:naade}.

\subsection{Biquaternions}\label{sec:biq}
The  biquaternions, introduced in \cite{clifford73},   consist of two copies of the quaternions,  with the eight units $\{1, \vec{i}, \vec{j}, \vec{k}, \epsilon, \epsilon \vec{i}, \epsilon \vec{j}, \epsilon \vec{k}\}$ where $\epsilon$ is a new unit commuting with everything and satisfying $\epsilon^{2} = 0$.\footnote{The cases $\epsilon^{2} = \pm 1$ were also considered by Clifford and lead to noneuclidean geometries (and their kinematics and rigid body mechanics), but lie outside the scope of this article. Eduard Study also made significant contributions to this field, under the name of \emph{dual} quaternions, which for similar reasons remain outside the scope of this article.} 

When one removes the constraint on the Euler top that its center of gravity remains fixed, one obtains the \emph{free top}, which, as the name implies, is free to move in space.  Here the allowable isometry group expands to the orientation-preserving euclidean motions $\Eucg{3}$, a six-dimensional group which is a semi-direct product of the rotation and translation subgroups.  Just as $SO(3)$ can be represented faithfully by the unit quaternions, $\Eucg{3}$ can be faithfully represented by the unit biquaternions.  

\subsubsection{Imaginary biquaternions and lines} Analogous to the way imaginary quaternions represent vectors in $\R{3}$, imaginary biquaternions represent arbitrary lines in euclidean space. An imaginary biquaternion of the form $a \vec{i} + b\vec{j}+c\vec{k}$ (a \emph{standard} vector) represents a line through the origin; one of the form $\epsilon(a \vec{i} + b\vec{j}+c\vec{k})$ (a \emph{dual} vector) represents an ideal line (aka \quot{line at infinity}).  An imaginary biquaternion whose standard and dual vectors are perpendicular (as ordinary vectors in $\R{3}$) corresponds to a line in $\RP{3}$; the general imaginary biquaternion represents a \emph{linear line complex},  the fundamental object of 3D kinematics and dynamics in this context. 

\subsubsection{Screw motions via biquaternion sandwiches} \label{sec:smvbs} For example, a general euclidean motion $\vec{R}$ is a screw motion with a unique invariant line, called its \emph{axis}.  $\vec{R}$  rotates around this axis by an angle $\alpha$ while translating along the axis through a distance $d$.  The axis can be represented by a unit imaginary biquaternion $\vec{m}$, as indicated above; then $\vec{v} := {\frac{\alpha+\epsilon d}{2}\vec{m}}$, a linear line complex, is the \emph{infinitesimal generator} of the screw motion and the exponential $e^{\frac{\alpha+\epsilon d}{2}\vec{m}}$ yields a unit biquaternion $\vec{g}$.  The sandwich operation $\vec{g}\vec{x}\vec{g}^{-1}$  gives the action of the screw motion on an arbitrary line  in space, represented by an imaginary biquaternion $\vec{x}$.  One can also provide a representation for its action on the points and planes of space: a point $(x,y,z)$ is represented as $1+\epsilon(x \vec{i}+y\vec{j}+z\vec{k})$ and a plane $ax+by+cz+d=0$ maps to $\epsilon d + a\vec{i} + b\vec{j} + c\vec{k}$, but then the sandwich operators have slight irregularities  (\cite{blaschke42}).  These irregularities are the necessary consequence of introducing \emph{ad hoc} representations for elements (points and planes) which do not naturally have a representation in the algebra.  We return to this point in \Sec{sec:lackcov} since it has generated some confusion related to our main theme.

\subsubsection{Seamless integration via ideal elements} \label{sec:biqideal} Note that the biquaternion representation seamlessly handles  cases which the traditional linear algebra approach has to handle separately.  For example, when the generating bivector $\vec{m}$ is an ideal line (a pure dual vector), then the isometry is a translation; hence one sometimes says, \quot{a translation is a rotation around a line at infinity};  in dynamics, a force couple (\emph{resp.}, angular momentum) is a force (\emph{resp.}, momentum) carried by an ideal line. Hence, one obtains the Euler top from the free top by constraining all momenta to be carried by ideal lines.  
  
\subsubsection{Euler equations of motion} The Euler equations of motion for the free top are then given by the same pair of ODE's given above, except that the symbols are to be interpreted in the biquaternion rather than quaternion context (where one uses unit (\emph{resp.}, imaginary) biquaternions instead of unit (\emph{resp.}, imaginary) quaternions).  These equations, as inherited by PGA, will be shown below  in \Sec{sec:krbdasct} to compare favorably with modern alternatives for kinematics and rigid body mechanics.
  
\subsubsection{Historical notes} The biquaternions, as developed  in \cite{study03}, were reformulated by von Mises (\cite{vm24}) using tensor and matrix methods; he called the result the \emph{motor algebra}.  The motor algebra has been developed further and applied in robotics by modern researchers, for example, \cite{corrochano00}.  \cite{ziegler} is an excellent historical monograph, recounting how the combined efforts of eminent mathematicians including M\"{o}bius, Pl\"{u}cker, Klein, Clifford, Study, and  others, led to the discovery of the biquaternions as a model for euclidean (and non-euclidean!)  kinematics and rigid body mechanics.  

\section{From exterior algebra to geometric algebras}
\label{sec:featga}
To understand how to  go beyond biquaternions to obtain a geometric algebra for euclidean geometry, we have need of  other mathematical innovations of the $19^{th}$ century: projective geometry and exterior algebra.  We first review some fundamental facts from projective geometry that are crucial to understanding the following treatment.  The remaining discussion in this section focuses on two  topics important to this exposition, not well-represented in the current literature: 
\begin{compactenum}[\textbullet]
\item the use of the \emph{dual exterior algebra} to construct geometric algebras not available otherwise, and 
\item the \emph{Cayley-Klein construction} of metric spaces atop projective space, particularly the delicate subject of degenerate signatures.
\end{compactenum}

\subsection{Preliminary remarks on vector space and projective space} First, recall that real projective space $\RP{n}$ can be derived from $\R{n+1}$ by introducing the equivalence relationship $\vec{x} \equiv \vec{y} \leftrightarrow \exists \lambda \neq 0$ with $\vec{x} = \lambda \vec{y}$, that is, the points of $\RP{n}$ are the lines through the origin of $\R{n+1}$; this construction is sometimes called the \quot{projectivization} of the vector space, and plays a large role in what follows.  We can either interpret an $n$-vector as being a vector in a vector space, or as representing a point in a projective space;  the former we refer to as the \emph{vector space} setting; the latter, as the \emph{projective} setting.
Also recall that to every real vector space $V$ there is associated the \emph{dual} vector space $V^{*}$, consisting of the linear functionals $V\rightarrow \mathbb{R}$.  The dual space in turn can be (and, in the following, is)  canonically identified with the  hyperplanes of $V$ by associating $f \in V^{*}$ with its  kernel $K_f := \{x \in V \mid f(x) = 0 \}$, a hyperplane. $\lambda f$ for $\lambda \neq 0$ is associated to the \emph{weighted} hyperplane $\lambda K_f$.

\subsection{Projectivized exterior algebra}\label{sec:projextalg} Begin with the standard real Grassmann or exterior algebra  $\grass{n+1}$ which encapsulates the subspace structure of the real vector space $\mathbb{R}^{n+1}$ (without inner product).  It is a graded associative algebra. 
The 1-vectors represent the vectors of $\mathbb{R}^{n+1}$.  The higher grades are constructed via the \emph{wedge} product, an anti-symmetric, associative product which is additive on the grade of its operands, and represents the \emph{join} operator on subspaces.  Projectivize this algebra  to obtain the projectivized exterior algebra $\pgrass{n+1}$ which in a natural way represents the subspace structure of real projective $n$-space $\RP{n}$, as built up out of points by joining them to form higher-dimensional subspaces.\footnote{One can also begin with begin with projective space and construct the Grassmann algebra in the obvious way; one obtains the same algebra $\pgrass{n+1}$.}

\subsection{Dual projectivized exterior algebra} The above process can also be carried out with the dual vector space $\RD{n+1}$ to produce the dual projectivized exterior algebra $\pdgrass{n+1}$.  It also models the subspace structure of $\RP{n}$,  but dually, so that the 1-vectors represent hyperplanes (by the canonical identity mentioned above between the dual space $\RD{n+1}$ and the hyperplanes of $\R{n+1}$).  The wedge product corresponds to the \emph{meet} or intersection of subspaces.  To avoid confusion we write the wedge operator in $\pdgrass{n+1}$ as $\wedge$ (meet) and the wedge operator in $\pgrass{n+1}$ as $\vee$ (join).  

\subsection{Duality} \label{sec:duality} For practical applications, it is necessary to be able to carry out both meet and join in a given exterior algebra. For example, consider the meet operator  in $\pgrass{n+1}$. This is often, in the context of a standard Grassmann algebra, called the \emph{regressive} product.\footnote{In a \emph{dual} Grassmann algebra, the regressive product is the \emph{join} operation. In general, it's the \quot{other} subspace operation not implemented by the wedge product of the algebra.}  To implement the meet operator  in $\pgrass{n+1}$,  we take advantage of Poincar\'{e} duality (\cite{greub2-67}, Sec. 6.8).  The basic idea is this: the exterior algebra and its dual provide two views on the \emph{same} projective space.  Any geometric entity in $\RP{n}$ appears once in each algebra. The Poincar\'{e} isomorphism is then a grade-reversing vector-space isomorphism $\vec{J}: \pgrass{n+1} \leftrightarrow \pdgrass{n+1}$  that maps a geometric entity of $\pgrass{n+1}$ to the \emph{same} geometric entity in $\pdgrass{n+1}$.  In this sense it is an \emph{identity} map; sometimes called the \emph{dual coordinate} map. Equipped with $\vec{J}$ we  define a meet operation $\wedge$ in $\pgrass{n+1}$ by \[ \vec{X} \wedge \vec{Y} := \vec{J}(\vec{J}(\vec{X}) \wedge \vec{J}(\vec{Y}))\] and similarly, a join operator $\vee$ for $\pdgrass{n+1}$.   An alternative non-metric method for calculating the join operator is provided by the shuffle product, see \cite{selig05}, Ch. 10. 

\subsection{Cayley-Klein construction of metric spaces} \label{sec:cayklngen}
Recall that the Sylvester Inertia Theorem asserts that a symmetric bilinear form $\vec{Q}$ can be characterised by an integer triple $(p,n,z)$, its \emph{signature}, describing the number of basis elements $\e{i}$ such that $\vec{Q}(\e{i}, \e{i}) = \{1,-1,0\}$, resp. 
If one attaches such a $\vec{Q}$ to $\RP{n}$ then,  for many choices of $\vec{Q}$, Cayley and Klein showed it is possible to define a distance function on a subset $\mathcal{M} \subset \RP{n}$ which makes $\mathcal{M}$ into a constant-curvature metric space (\cite{klein:neg}, Ch. 6, \cite{gunn2011}, \S 3.1, or \gTh, Ch. 4). When $z=0$, the construction is based on the projective invariance of the cross ratio of four collinear points.  The distance between two points $\vec{A}$ and $\vec{B}$ is defined as \[d(\vec{P}, \vec{Q}) := k\ln({CR(\vec{A}, \vec{B}; \vec{F}_1, \vec{F}_2)})\]  where $k$ is an appropriate real (or complex) constant, and $\vec{F}_1$ and $\vec{F}_2$ are the two intersection points (real or imaginary) of their joining line with the absolute quadric (the null vectors of $Q$).  Since the cross ratio is a multiplicative function, $d$ is additive, and satisfies the other properties of a distance function.

For example the signature $(n+1,0,0)$ leads to elliptic space $\Ell{n}$, while $(n,1,0)$ leads to hyperbolic space $\Hyp{n}$.  \textbf{Note:} in the vector space setting,  the signature $(n+1,0,0)$ produces  the \emph{euclidean} metric; in the projective setting (via the Cayley-Klein construction), however, the same signature produces the \emph{elliptic} (or \emph{spherical}\footnote{Two copies of elliptic space $\Ell{n}$ can be glued together to obtain spherical space $\Sph{n}$.}) metric.  This ambiguity has led to  misunderstandings in the literature which we discuss below in \Sec{sec:wrongmetric}.

\subsubsection{Cayley-Klein construction of euclidean space} 
\label{sec:caykln}
The signature for euclidean geometry is degenerate, that is, $z\neq0$.  Since this is a crucial point, we motivate the correct choice using the example of the euclidean plane  $\Euc{2}$. Consider two  lines $m_{1} : a_{1}x+b_{1}y+c_{1}=0$ and $m_{2} : a_{2}x+b_{2}y+c_{2}=0$.  Assuming WLOG $a_{i}^{2}+b_{i}^{2} = 1$, then $\cos{\alpha} = a_{1}a_{2}+b_{1}b_{2}$, where $\alpha$ is the angle between the two lines: changing the $c$ coefficient translates the line but does not change the angle it makes to other lines. This generalizes to the angle between two hyperplanes in $\Euc{n}$.  One coordinate  plays no role in the angle calculation, hence the signature has $z=1$.  Thus, the Cayley-Klein construction applies the signature $(n,0,1)$ to \emph{dual} projective space to obtain a model for euclidean geometry in $n$ dimensions.  

\subsubsection{Euclidean distance between points}  \label{sec:eucdist} The discussion above takes its starting point so that the angle between lines (hyperplanes) can be calculated.  What about the distance between points?  Already Klein (\cite{klein:neg}, Ch. 4 \S 3) provided an answer to this question (or in English see \gTh, \S 4.3.1).  The inner product given above on lines (hyperplanes) induces the (very degenerate) signature $(1,0,n)$ on points, so that one cannot measure the distance between points via the inner product.  However, using a sequence of non-euclidean signatures that converge to $(n,0,1)$ (\gTh, \S 3.2), one can show that the distance function between normalized homogeneous points $\vec{P}$ and $\vec{Q}$ converges, up to an arbitrary positive constant\footnote{This positive constant determines the length scale of euclidean space, such as feet or meters.}, to the familiar euclidean formula $\| \vec{P} -\vec{Q} \|$, the length of the euclidean vector $\vec{P} - \vec{Q}$ (here we apply \emph{euclidean} in the vector space setting).  

\subsection{Geometric algebras from Cayley-Klein} 

It is straightforward to obtain geometric algebras from the Cayley-Klein construction.\footnote{In fact, given the personal and scientific  friendship of Klein and Clifford in the 1870's, it is likely that the Cayley-Klein construction influenced both Clifford's discovery of biquaternions (1873) and geometric algebra (1878) (\cite{ziegler}, Ch. 7).}
We use the signature of the Cayley-Klein construction to define the inner product for our geometric algebra.  To remind us that we are  operating within the projective space setting rather than the vector space one, we call such a GA a \emph{projective geometric algebra} or PGA for short. The name reflects the fact that in such a geometric algebra,  the metric is based (either directly or via a limiting process) on the \emph{projective cross ratio} (as explained in \S \ref{sec:cayklngen} above).  
Hence, a projective geometric algebra is a special case of a homogeneous geometric algebra (\cite{dfm07}, Ch. 11), which is also sometimes called a \quot{1-up} model, since it requires $(n+1)$-dimensional coordinates to represent the an $n$-dimensional metric space.   
So, $\pclal{n+1}{1}{0}$ (\emph{resp.}, $\pclal{n}{1}{0}$) provides a model for $n$-dimensional elliptic (\emph{resp.}, hyperbolic) space, and is called elliptic (\emph{resp.}, hyperbolic) PGA. 

\subsubsection{Dual geometric algebras} \label{sec:dga} It is also possible to use the dual exterior algebra $\pdgrass{n+1}$ as the basis for a geometric algebra.  Thus, the inner product is defined on the hyperplanes of the projective space (the 1-vectors).   We call such a geometric algebra a \emph{dual} geometric algebra;  a geometric algebra built atop $\pgrass{n+1}$ we call \emph{standard} to distinguish it from the dual case. We sometimes call the former a \emph{plane-based} algebra and the latter a \emph{point-based} one, emphasising the very different meaning of the 1-vectors in the respective cases. One can compare the standard and dual GA with the same signature by calculating the induced metric on $n$-vectors (which correspond to the 1-vectors in the dual algebra). One finds that the dual algebra $\pdclal{n+1}{0}{0}$ yields elliptic space again. 
$\pdclal{n}{1}{0}$, on the other hand, yields dual hyperbolic space, built up of the hyperplanes lying \emph{outside} the unit sphere (rather than points \emph{inside} the unit sphere). The induced signature $(n,1,0)$ on $n$-vectors  (calculated by writing the basis $n$-vectors as products of 1-vectors and squaring the results) is, however, the same as in the standard algebra $\pclal{n}{1}{0}$.  When the metric is non-degenerate, as here,  the dual geometric algebra can  be obtained by multiplying the original algebra by the pseudoscalar $\eye$ and then reversing the grades. 
Hence, a non-degenerate signature applied to the dual exterior algebra yields nothing new; every metric relationship in the dual algebra is mirrored in the standard algebra via pseudoscalar multiplication. 
This is not true for degenerate signatures, as we see in the next section.

\subsection{The degenerate signature $(n,0,1)$} \label{sec:eucspace}  We established above that the degenerate signature $(n,0,1)$ applied to the dual Grassmann algebra leads to the euclidean algebra $\pdclal{n}{0}{1}$. The standard algebra $\pclal{n}{0}{1}$, however, represents a different metric space, \emph{dual euclidean space}. These \emph{cannot} be obtained from one another by pseudoscalar multiplication since the pseudoscalar is not invertible. For example, for two normalized $k$-blades $\vec{A}$ and $\vec{B}$, $\vec{A}$ and $\vec{B}$ are parallel $\leftrightarrow \vec{A}\eye = \vec{B} \eye$.  The induced signature on $n$-vectors, $(1,0,n)$, is very degenerate, and not equivalent to the signature on 1-vectors. As a result, euclidean and dual euclidean space exhibit an asymmetry not present in the non-degenerate case: the absolute quadric of  euclidean space is a single ideal \emph{plane}, while that of dual euclidean space is a single ideal \emph{point}.  This reflects the fact that  euclidean space arises by letting the curvature of a non-euclidean space go to $0$, while dual euclidean space arises when the curvature goes to $\infty$. 

\subsubsection{Dual euclidean space}  \label{sec:des}  
Because the distinction between euclidean and dual euclidean space is crucial to the theme of this article, and is not well-known,  we discuss it briefly here.  The simplest example of a dual euclidean space occurs  within the hyperbolic algebra $\pclal{n}{1}{0}$ (which forms the basis of conformal geometric algebra, \S \ref{sec:cga} below).  A hyperplane tangent to the null sphere $\vec{Q}$ at a point $\vec{P}$ has induced signature $(n,0,1)$.  $\vec{P}$ provides the degenerate basis vector  satisfying $\vec{P}^{2}=0$, all other points have non-zero square since they do not lie on $\vec{Q}$.  
Furthermore, no standard geometric algebra can contain euclidean space as a flat subspace in this way.  Why?  We saw above that the induced signature $(1,0,n)$ on points is more degenerate than the signature $(n,0,1)$ on hyperplanes.  This asymmetry is incompatible with an algebra in which 1-vectors represent points;  only a \emph{dual} geometric algebra can provide both the required signature $(n,0,1)$ on hyperplanes and $(1,0,n)$ on points.  
Dual euclidean space shows promise as a tool for effectively modeling some aspects of the natural world, see \cite{kowol09} and \gTh, Ch. 10.

\subsection{Geometric algebras for euclidean geometry}
 In this section we give an overview of the field of candidates of geometric algebras for doing euclidean geometry.  We have already met one of the candidates, $\pdclal{n}{0}{1}$, above in \Sec{sec:eucspace}; we describe in it more detail now.  
 We will see below in \ref{sec:wrongmetric} that homogeneous models with non-degenerate metrics are inferior to $\pdclal{n}{0}{1}$ for euclidean geometry.  The other remaining candidate for \thetask is a 2-up model,  conformal geometric algebra (CGA), which we introduce next.   
  
\subsubsection{Projective geometric algebra}\label{sec:pga}

The geometric algebra  $\pdclal{n}{0}{1}$ introduced above  for euclidean geometry we call \emph{euclidean} PGA.  When the context makes it clear, as generally in the remainder of this article, we refer to it simply as PGA.    Other examples of PGA's are elliptic PGA ($\pclal{n+1}{0}{0}$) and hyperbolic PGA ($\pclal{n}{1}{0}$).

The measurement of angles is given then by the inner product on the 1-vectors as described above in \ref{sec:caykln}.  The distance function between points, also described there, appears in (at least)  two different sub-products of the algebra: $d(\vec{P}, \vec{Q}) = \| \vec{P} \vee \vec{Q} \| = \| \vec{P} \times \vec{Q} \|$ (assuming that $\vec{P}$ and $\vec{Q}$ have been normalized).   Here, $ \vec{P} \vee \vec{Q}$ is the joining line of the points, and $ \vec{P} \times \vec{Q} := \langle \vec{P} \vec{Q}\rangle_{2}$ is the orthogonal complement of the joining line.  Details of the first of these formulas, and many other formulas, can be found in \gTh, Ch. 6 and Ch. 7.  The absolute quadric is the \emph{ideal plane};  because of its importance we introduce for it the notation $\vec{\omega}_{P}$ ('P' stands for \emph{projective}).

For $n=3$, the case of most general interest, the even subalgebra $\pdclplus{3}{0}{1}$ is isomorphic to the biquaternions.  To construct the isomorphism, map the imaginary biquaternions to the bivectors of $\pdclplus{3}{0}{1}$ in the obvious way (since both provide Pl\"{u}cker coordinates for line space), and $\epsilon$ to the pseudo-scalar $-\eye$ of the geometric algebra. 
This isomorphism brings with it the elegant  representation of rigid body motion described above in \Sec{sec:biq}.  The representation can be extended to include points and planes; details can be found in \cite{gunn2011}, \S 15.6.  Also note that PGA replaces the irregular transformation formula for the sandwich operators of the biquaternions and  of the motor algebra (\Sec{sec:smvbs}), with the uniform sandwich operators of the geometric algebra. \textbf{Warning:} The biquaternions are also isomorphic to $\pclplus{3}{0}{1}$, the even subalgebra of dual euclidean space.  We return to this later in the article (\Sec{sec:lackcov}) as it appears to have been a source of confusion relevant to our theme. 

$\pdclal{3}{0}{1}$ has the distinction of integrating two of William Clifford's most significant inventions, geometric algebra and biquaternions, into a single algebra.\footnote{That Clifford himself appears to have overlooked this algebra is not surprising, considering the tentative nature of his research into both of these objects (unavoidable due to his early death);  the presence of a degenerate signature and the use of a dual exterior algebra are both features of geometric algebra which were not known during his lifetime.}
Seen in this light, $\pdclal{3}{0}{1}$ stands in the confluence of two streams of $19^{th}$ century mathematics:  on the one hand, that leading to the \MN biquaternion formulation of rigid body mechanics, and on the other hand, the Cayley-Klein integration of metric geometry in projective geometry, so that it has close connections to the genesis of geometric algebra itself. 
The algebra first appeared in the modern literature in \cite{selig00} and \cite{selig05}, and was then extended and embedded in the \MN toolkit described in \gTh.  A compact, self-contained treatment is given in \cite{gunn2011} (extended version \cite{gunnFull2010}).  

\subsubsection{Conformal geometric algebra} \label{sec:cga}
If one begins with the $(n+1)$-dimensional PGA $\pclal{n+1}{1}{0}$ for hyperbolic geometry, one can obtain another model for euclidean geometry as follows.  Identify the points of the absolute quadric $\vec{Q}$ (the null sphere) with  $\Euc{n}$ by stereographic projection.  Then one can normalize the coordinates of these points so that the inner product between two points yields the square of the euclidean distance between the two points.  Points outside (inside) the sphere can be identified with spheres in $\R{n}$ of positive (negative) radius. The points of the null sphere itself can be identified with the points of $\Euc{n}$ itself; and are sometimes called \emph{zero-radius spheres}.   Projectivities which preserve  $\vec{Q}$ correspond to conformal maps of $\R{n}$, hence this model is called the \emph{conformal} model of euclidean geometry, and the associated geometric algebra is called conformal geometric algebra (CGA).  It was introduced in its present form in \cite{hlr01} and has developed rapidly since then (\cite{doran03}, \cite{dfm07}, \cite{perwass09}, and \cite{ggap2011}) . In light of the prolific literature available, we omit a more detailed description here.  

\myboldhead{The flat representation in CGA} CGA contains a sub-algebra closely related to PGA.  Since it will play a role in the sequel we describe it here. As noted above, the tangent plane $T_\vec{P}$ at a point $\vec{P}$ of the null sphere of CGA is a sub-algebra isometric to dual euclidean space $\pclal{n}{0}{1}$.  Letting $\vec{P} = n_\infty$, and polarizing $T_\vec{P}$ by multiplying by the pseudoscalar $\eye \in \pclal{n+1}{1}{0}$ yields another sub-algebra isometric to euclidean space $\pdclal{n}{0}{1}$.  We call this sub-algebra $\mathfrak{S}$.  It consists of all flat subspaces containing $n_\infty$. The relationship to PGA is this: in PGA, the $k$-dimensional subspaces of $\Euc{n}$ are represented by $(n-k)$-vectors of the algebra; in $\mathfrak{S}$, the $k$-dimensional subspaces of $\Euc{n}$ are represented by the $(k+1)$-vectors containing the \quot{star} point $n_\infty$.  The name \emph{flat} representation comes from the fact that one can also obtain it by taking the standard representation of CGA (as zero-radius spheres) and wedge it with $n_\infty$.  For the purposes of this article, we content ourselves with the observation that $\mathfrak{S}$ and $\pdclal{n}{0}{1}$ are isometric, hence are essentially identical (except that $\mathfrak{S}$ has an extra, irrelevant dimension). Further work to establish the exact relationship between these two representations needs to be done.  


\section{Clarification work}
\label{sec:clar}
We have identified the two algebras PGA and CGA as candidates \thetask.  There exists some controversy in the literature whether there might be other candidates, as well as questions regarding the suitability of PGA.  We now turn to examine these issues in more detail.   In this section we dispose of other  homogeneous models which appear in the literature.  

\subsection{Which homogeneous model?}
\label{sec:wrongmetric}
Here we want to discuss other homogeneous (i. e., 1-up) models for the chosen task besides $\pdclal{n}{0}{1}$.   
For example, Chapter 11 of \cite{dfm07}, entitled \emph{The homogeneous model}, describes one such model 
(which appears in several other  textbooks (\cite{doran03}, \cite{perwass09}). 
In \S 11.1 one reads:
\begin{quote}
...[the homogeneous model of euclidean geometry] embeds $\R{n}$ in a space $\R{n+1}$ with one more dimension and then uses the algebra of $\R{n+1}$ to represent those elements of $\R{n}$ in a structured manner.
\end{quote}
The authors then reject the use of a degenerate metric as \quot{ inconvenient},  and therefore propose using any non-degenerate metric, for example,  $(n+1,0,0)$ or $(n,1,0)$, yielding the geometric algebras $\pclal{n+1}{0}{0}$ and $\pclal{n}{1}{0}$.  That means the basis element $\vec{e}_{0}$ satisfies $\vec{e}_{0}^{2} = \pm 1$. 
%
We mentioned above in \Sec{sec:cayklngen} that these two algebras yield a \emph{elliptic}, resp. \emph{hyperbolic},  metric on projective space.  Here we apparently encounter the widespread confusion between the meaning of \quot{euclidean} in the vector space versus the projective space setting.  We return to the consequences of this confusion below in \Sec{sec:rephrase}. 

\subsubsection{Comparison based on worked-out example}\label{sec:example}
\begin{figure}
 \def\xyz{.31}
 \def\zyx{.01}
{\setlength\fboxsep{0pt}\fbox{\includegraphics[width=\xyz\columnwidth]{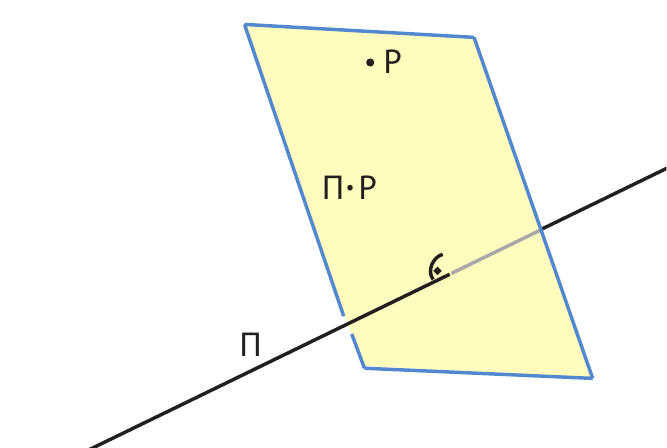}}}\hspace{\zyx\columnwidth}
{\setlength\fboxsep{0pt}\fbox{\includegraphics[width=\xyz\columnwidth]{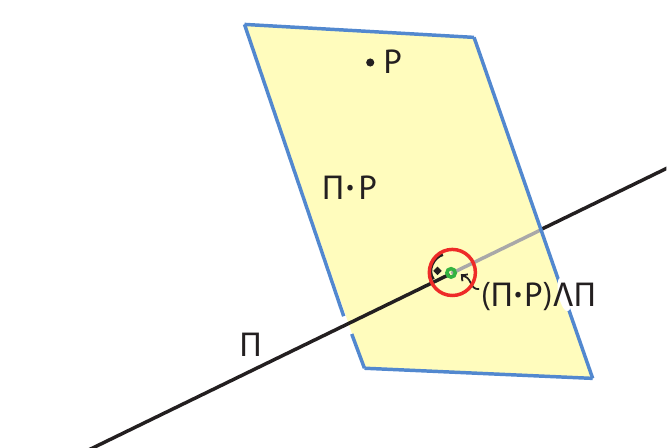}}}\hspace{\zyx\columnwidth}
{\setlength\fboxsep{0pt}\fbox{\includegraphics[width=\xyz\columnwidth]{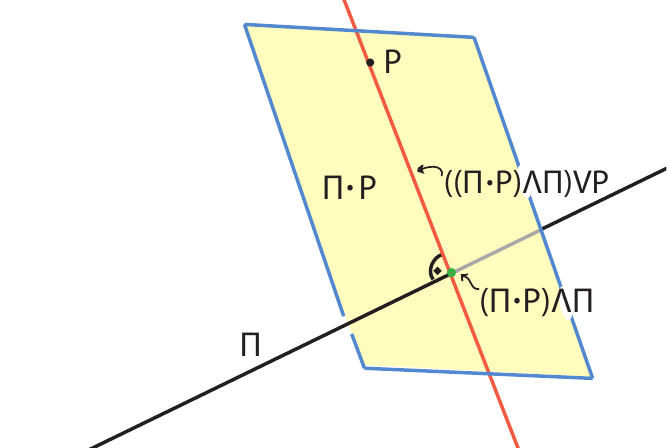}}}
\vspace{-.2in}
\caption{Solution of the exercise in PGA.} 
\label{fig:perpPtLn}
\end{figure}

We  compare this non-degenerate homogeneous model with $\pdclal{n}{0}{1}$ on a simple geometric construction,  taken from \S 11.9 of \cite{dfm07}:
\begin{quote}
Given a point $\vec{P}$ and a non-incident line $\momo$ in $\Euc{3}$, find the unique line $\sigo$ passing through $\vec{P}$ which meets $\momo$ orthogonally.
\end{quote}
$\pdclal{3}{0}{1}$ yields directly the compact solution $((\momo \cdot \vec{P} )\wedge \momo)  \vee \vec{P}$;    \Fig{fig:perpPtLn} decomposes the solution in three easy-to-understand steps.  This PGA solution is coordinate-free and \MN,  hence valid for hyperbolic and elliptic space also. The first step is the most important: $\momo \cdot \vec{P} = \langle \momo \vec{P} \rangle_{1}$ is the plane $\vec{p}$ through $\vec{P}$ perpendicular to $\momo$.
 
 We adopt the solution from \cite{dfm07}, p. 310,  to conform to the notation used here (so $\vee$ is join, $\wedge$ is meet, $\cdot$ is contraction, $\vec{A}^{\perp} = \vec{A} \eye$ is the polar of $\vec{A}$).  The reasoning is similar.  Once the perpendicular plane $\vec{p}$ has been produced, the desired line $\sigo$  can be obtained, as in PGA, using $(\vec{p} \wedge \momo) \vee \vec{P}$.  The difference lies in the definition \[\vec{p} :=  (\vec{P} \cdot ((\e{0} \cdot \momo) \vee \e{0}))^{\perp}\]  (We have chosen $\e{0}^{2}=1$ to simplify the expressions.)  The reader can verify that $(\e{0} \cdot \momo) \vee \e{0}$ is the line through the origin $\e{0}$ parallel to $\momo$.  Hence when $\momo$ passes through $\e{0}$, the expression obtained is the same as that in PGA (modulo the presence of the polar operator ${ }^{\perp}$, which reflects the fact that we are working in a standard rather than dual GA). When $\momo$ doesn't pass through $\e{0}$, then by translating $\momo$ there, one obtains the desired answer, since a plane through $\vec{P}$ perpendicular to a line through $\e{0}$ will also be perpendicular to any translated copy of this line. 
 
 Here we see why $\e{0}$ appears in \emph{every} expression obtained in the discussion in \cite{dfm07}
 : exactly at $\e{0}$, the elliptic metric and the euclidean metric agree.  This is equivalent to the fact, that only at $\e{0}$ does the polar plane  in the elliptic metric agree with the euclidean polar plane (which is always the ideal plane).  This  allows one to translate geometric entities to the origin, operate on them in a \MN way (for example to obtain their euclidean directions), and translate them back when necessary. In comparison to the PGA solution, however, the one obtained in this way is neither coordinate-free nor \MN -- in addition to  involving an extra pair of operations to translate the line to the origin.   

Things become even less satisfactory when one attempts to implement euclidean isometries.   The authors acknowledge that using the non-degenerate metric,  it is impossible to express euclidean translations as sandwich operators (\S 11.8 of \cite{dfm07}).
This leads to the conclusion:
\begin{quote}
The main problem with using the metric of $\R{n+1}$ is that you cannot use it directly to do Euclidean geometry, for it has no clear Euclidean interpretation.
\end{quote}
The foregoing quote is a good motivation for  the next section, where we attempt to clarify the situation by differentiating various meanings of $\R{n}$ and \quot{euclidean}. 
 
\subsubsection{Three meanings of $\R{n}$} 
\label{sec:threern}
The symbol $\R{n}$ occurs 5 times in the two short quotes of the above section, grounds for asking what exactly it means.
We can in fact distinguish at least three different meanings:
\begin{compactenum} [1)]
\item \hspace{-.1in} \textbf{Vector space.}  In this form,  $\R{n}$ represents the vector space used to define real projective space $\RP{n-1}$. 
It is an $n$-dimensional linear space  with an addition operation, real scalar multiplication, and distributive law, but \textbf{without} inner product. One can develop a theory of linear mappings between such spaces, and from this,  the dual vector space $\VSD$. The \emph{evaluation map} $\R{n} \otimes \RD{n} \rightarrow \mathbb{R}$ of a vector and a dual vector (linear functional), often written $\langle \vec{v}, \vec{\mu} \rangle := \vec{\mu}(\vec{v})$ is sometimes confused with an inner product. 
We recommend using the terminology ($n$-dimensional) $\VS$ for this meaning of $\R{n}$, whenever possible.\footnote{It is not always convenient, see for example \ref{sec:projextalg} above, where $\R{n}$ is traditionally used to define the exterior algebra, even though the inner product plays no role thereby.}   See \cite{greub1-67}, Chapter 1-2 for details.
\item \hspace{-.1in} \textbf{Inner product space.} One begins with a vector space and adds an \emph{inner product } between pairs of vectors, which is a  symmetric bilinear form on the vector space. This produces an \emph{inner product space}.  When the form is positive definite, it's called a   \emph{euclidean} inner product space. We recommend retaining the use of  $\R{n}$ for this meaning.  Consult  \cite{greub1-67}, Chapter 7 for details on inner product spaces.
\item \hspace{-.1in}  \textbf{Euclidean space.} This is a simply-connected  \emph{metric space}, of constant curvature 0, homeomorphic to $\R{n}$ but equipped with the Euclidean distance function (discussed for example in \gTh, Chapter 4) 
between its points.   We recommend using the notation $\Euc{n}$ for this  space.  The points of $\Euc{n}$ are in a 1:1 correspondence to the \emph{vectors} of $\R{n}$ (the origin of $\Euc{n}$ maps to the zero vector of $\R{n}$), but  $\Euc{n}$ is \textbf{not}  a vector space, and the inner product discussed in the previous item has, \emph{a priori}, nothing to do with the measurement of distances in $\Euc{n}$.  
\end{compactenum} {

Armed with these three different meanings which sometimes are attached to the same symbol $\R{n}$, along with two meanings of \emph{euclidean} (depending on whether one is in the \emph{vector space} or the \emph{projective} setting), let's return to the discussion of the homogeneous model.
%
\subsubsection{Rephrasing using the differentiated notation}\label{sec:rephrase}
When we apply this differentiated terminology and what we learned about the Cayley-Klein construction of metric spaces to the initial quote from \cite{dfm07}, we arrive at the following:
\begin{quote}
...[the homogeneous model] embeds $\Euc{n}$ in a real vector space $\VS$ of dimension $n+1$ and then uses the algebra of $\VS$ to represent the elements of $\Euc{n}$ in a structured manner.
\end{quote}
In this form, $\R{n}$ no longer occurs: there is no longer a given real vector space nor inner product, implied by the original definition.  
Consequently, one can use this  modified description as a {starting} point for the \emph{search} for the correct choice of Cayley-Klein space; we have  sketched  above how one arrives at the dual vector space $\VS = \RD{n+1}$ with signature $(n,0,1)$, yielding the algebra $\pdclal{n}{0}{1}$. 

To sum up: this confusion of the three meanings of $\R{n}$ and two meanings of \emph{euclidean} means that many of the objections to \quot{the} homogeneous model appear in the light of the foregoing discussion as legitimate complaints against using the \emph{wrong vector space} or the  \emph{wrong signature} to model euclidean geometry. In the next section we turn to consider if there are other choices which yield better results.

\section{Homogeneous models using a degenerate metric}
\label{sec:hmudm}
Faced with the difficulties ensuing on the use of the non-degenerate metric, \cite{dfm07}, p. 314, states:
\begin{quote}
We emphasize that the problem is not geometric algebra itself, but the homogeneous model and our desire to use it for euclidean geometry. It will be replaced by a much better model for that purpose in Chapter 13 [the conformal model - \emph{cg}].
\end{quote}

In fact, what the authors of  \cite{dfm07} have shown is that a \emph{homogeneous model with non-degenerate metric} is \quot{the problem} -- recall that  the use of degenerate metrics was rejected as inconvenient. Hence it remains to be seen whether a PGA based on a degenerate metric, such as $\pdclal{n}{0}{1}$, could provide a faithful model for euclidean geometry.  We now turn to an analysis of three common objections to the use of such a degenerate metric.

\subsection{Objection 1: lack of covariance}
\label{sec:lackcov}
\emph{Covariance} has a variety of meanings related to the behavior of maps and coordinate systems;  in our context, it is equivalent to the existence of sandwich (or \emph{versor}) implementations of the euclidean group $\Eucg{n}$ (\cite{dfm07}, p. 369). That such versors exist for the even subalgebra of $\pdclal{3}{0}{1}$ (and the associated Spin group) follows from its isomorphism with the biquaternions mentioned above in \Sec{sec:pga}.  
The extension to the full algebra (and the associated Pin group) is straightforward and is described in \cite{gunn2011}, \S 15.4.2, \S 15.5.4.


One can gain a different impression, however, from some of the current literature.  
For example, \cite{li08}, p. 11, identifies  the Clifford algebra $\pclal{n}{0}{1}$ (in our notation) as the appropriate algebra for euclidean geometry. It focuses on the case of $n=3$ and remarks on the isomorphism of the dual quaternions with the even subalgebra.  This leads to the remark:
\begin{quote}
However, the dual quaternion representations of primitive geometric objects such as points, lines, and planes in space are not \emph{covariant}. More accurately, the representations are not tensors, they depend upon the position of the origin of the coordinate system irregularly. 
\end{quote}

In the first place, the proper algebra for euclidean geometry is not  $\pclal{n}{0}{1}$ but the dual version $\pdclal{n}{0}{1}$ (see above, \Sec{sec:des}).  
This confusion is perhaps due to the fact that the dual quaternions are isomorphic to the even sub-algebra of \emph{both} these geometric algebras (\Sec{sec:pga} above).
Furthermore, as discussed above in \Sec{sec:smvbs}, the dual quaternion representation of points and planes is \emph{not} the same as the representation of points and planes in these geometric algebras: the dual quaternion representation has irregularities not exhibited by the geometric algebra representation due to the fact there there is no natural representation for points and planes within the algebra.  These irregularities however have to do with the form  the sandwich operators take and do not effect the covariance of the representation.  We are aware of no grounds  for the claim made here that the dual quaternion representations are not covariant.  It might  be due to mixing up the two algebras $\pclal{n}{0}{1}$ and $\pdclal{n}{0}{1}$ in the calculation, since each provides, taken for itself, covariant representations for their respective isometry group.

Readers who would like to confirm the claimed covariance for themselves, and do not have access to \cite{gunn2011} or \gTh are referred to Appendix A, which provides a detailed  discussion of the versors of $\pdclal{2}{0}{1}$ and their associated sandwich operators.

\subsection{Objection 2: lack of \emph{duality}} \label{sec:compdual}

Another common objection to the use of degenerate metrics is often expressed in terms of a \quot{lack of duality}. Consider the following quote from \cite{hlr01}, an article often  associated with the birth of modern CGA (\cite{dfm07}, \S 13.8):
\begin{quote}
Any degenerate algebra can be embedded in a non-degenerate algebra of larger dimension, and it is almost always a good idea to do so.  Otherwise, there will be subspaces without a complete basis of dual vectors, which will complicate algebraic manipulations.
\end{quote}
As with the case of \emph{euclidean} above, there are multiple meanings  for the term \emph{dual} in the literature which must be carefully differentiated.  Here there are at least two distinct meanings:
\begin{enumerate}
\item  The ability to calculate the \emph{regressive product}.  As shown in \Sec{sec:duality}, for this operation one only requires the \emph{dual coordinates} of a geometric entity, and this is provided in a non-metric way by  Poincar\'{e} duality in the exterior algebra, or by the shuffle product. 
\item The effect of multiplication by the pseudoscalar: $\Pi: \Pi(\vec{X}) := \vec{X} \vec{I}$.  The action of $\Pi$ on the  underlying projective space is known in classical projective geometry as the \emph{polarity on the metric quadric}, that is, a correlation (maps points to hyperplanes and vice-versa) that maps a point to its orthogonal  hyperplane with respect to the inner product encoded in the quadric, and vice-versa.\footnote{If one considers the inner product as a symmetric bilinear form $B(\vec{u}, \vec{v})$, then one obtains a linear functional ${B}_{u}$ by fixing $\vec{u}$ and defining  ${B}_{u}(\vec{v}) := B(\vec{u},\vec{v})$. Then the kernel of $B_{u}$ is the indicated orthogonal hyperplane: the set of all vectors with vanishing inner product with $\vec{u}$. This is also sometimes referred to in GA as the \emph{inner product null space} (IPNS) of $\vec{u}$.} When the metric is non-degenerate, then $\eye^{2} = \pm1$, and the polarity is a grade-reversing algebra bijection which is a vector-space isomorphism on each grade, and whose square is the identity\footnote{Since even if $\Pi^2=-1$,  $\vec{X}$ and $-\vec{X}$ represent, projectively, the same element.}. 
In this case, one can define the regressive product (analogously to the use of Poincar\'e duality in \Sec{sec:duality}): \[\vec{X} \wedge \vec{Y} := \Pi(\Pi(\vec{X}) \vee \Pi(\vec{Y}))\] 
\end{enumerate}

We recommend  that the term \emph{polarity} be adapted also in  geometric algebra  for pseudoscalar multiplication, to distinguish it from the previous non-metric meaning of \emph{duality}.  So that, in the above quote, one would speak of the \emph{polar} basis instead of the \emph{dual} basis.  The term \emph{dual} basis would be reserved for the result of the dual coordinate map $\vec{J}$. 

The implicit use of the metric to calculate the regressive product has a long tradition, going all the way back to Grassmann and continuing up to  \cite{zieghest91}, an influential modern article devoted to doing projective geometry using geometric algebra; its  continued use of $\Pi$ for the regressive product -- despite the absence of a natural metric for projective geometry -- appears to have cemented the misunderstanding described here. 
As a result, in the above quote as well as other popular texts (\cite{dfm07}, \cite{perwass09}, \cite{doran03}) one might falsely gain the impression that a non-degenerate metric  \emph{must} be used to implement the regressive product, as few or no details of an alternative are provided. 
But in the absence of an invertible pseudo-scalar, one always has access to Poincar\'{e} duality.  Hence this objection to PGA cannot be sustained.

\subsection{Objection 3:  Absence of invertible pseudoscalar}
\label{sec:absinv}
\cite{li08}, also p. 11, raises a further (related) objection to the use of a degenerate metric:
\begin{quote}
Because the inner product in $\mathbb{R}_{n,0,1}$ is degenerate, many important invertibilities in non-degenerate Clifford algebras are lost.\footnote{Since it is irrelevant to this object, we overlook the fact, discussed above, that the correct Clifford algebra here should be $\mathbb{R}_{n,0,1}^{*}$.}
\end{quote}
  
We have discussed the invertibility of $\eye$ as a condition for duality above in Objection 2, and shown that there are other means to implement duality.
It is true that many  formulas in the GA literature tend to be given in terms of $\eye^{-1}$ rather than $\eye$.  But in the cases we are familiar with, it is also possible to use $\eye$ instead, perhaps at the cost of a more complicated expression for the sign.   For example, \cite{hestenes10} defines the \emph{dual} (our \emph{polar}!) of a multivector $A^{*} := A\eye^{-1}$ and shows then that the inner and outer products obey the relation: $a\cdot A^{*} = (a\wedge A)^{*}$.  If instead one defines $A^{*} := A\eye$, the stated relation remains true, and valid for \emph{any} pseudoscalar.

Experience leads us, in fact, to a very different view of the non-invertible pseudoscalar: it has proved to be an \emph{advantage}, since it faithfully represents the metric relationships within euclidean geometry. The  calculation from \Sec{sec:example} provides a good example.  Consider the sub-expression $\momo \cdot \vec{P}$. Letting the point $\vec{P}$ move freely, one obtains a set of parallel planes $\momo \cdot \vec{P}$ which all have the same polar point $(\momo \cdot \vec{P})\eye$,  the ideal point of the line $\momo$.   For invertible pseudo-scalars,  however,  the polar points of distinct planes are distinct.  Or, recall the discussion of the special  role of $\e{0}$ in the solution  in \Sec{sec:example}: it does not appear in the simpler PGA formula, since in PGA, $\vec{P}\eye = \vec{\omega}_P$ for all normalized $\vec{P}$, while in the elliptic metric this is only true for $\vec{P} = \e{0}$.  

\section{Comparison: A feature-set for \quot{doing geometry}}
\label{sec:comp}
Having sketched its mathematical lineage dating back to Klein and Clifford, and disposed of a series of modern  objections which have been raised against it, the reader is hopefully convinced that euclidean PGA deserves the title of \quot{standard} or \quot{classical} homogeneous model of euclidean geometry.
We are now prepared to compare it to  the conformal model, CGA.  As a basis for this comparison, we rely on a recent tutorial on the conformal model \cite{dorst2011}. This tutorial describes  the challenge of \quot{doing geometry} on a computer, a challenge which matches well with \thetask we set out at the beginning of the article, so we use this tutorial as a basis for a first comparison of the two algebras. 
 The tutorial lists seven \quot{tricks} and three \quot{bonuses} which the conformal model offers in this regard.  We list them here:
\begin{compactenum}
\item Trick 1: Representing euclidean points in Minkowski space.
\item Trick 2: Orthogonal transformations as multiple reflections in a sandwiching representation.
\item Trick 3: Constructing elements by anti-symmetry.
\item Trick 4: Dual specifications of elements permits intersection.
\item Bonus: The elements of euclidean geometry as blades.
\item Bonus: Rigid body motions through sandwiching.
\item Bonus: Structure preservation and the transfer principle.
\item Trick 5: Exponential representation of versors.
\item Trick 6: Geometric calculus.
\item Trick 7: Sparse implementation at compiler level.
\end{compactenum}

How does  $\pdclal{n}{0}{1}$ stand with respect to these features?
In fact, it offers offers all the ten features listed.  Some slight editing is required to \quot{translate} to PGA; for example, Trick 1 has to be rephrased as \quot{Representing euclidean points in projective space}. Duality (trick 4) is implemented in a non-metric way in our homogeneous model, and is used to represent join, not intersection. There are naturally some elements of euclidean geometry which cannot be represented as blades in PGA (bonus 1), such as point pairs and spheres.  But the basic flat elements belonging to classical euclidean geometry are present: points, lines, and planes; and these are the ones belonging to  \thetask.  We return to the richer class of primitives in CGA in  \Sec{sec:round} below. 


One immediate corollary is that euclidean PGA  ($\pdclal{3}{0}{1}$) is the \emph{smallest known algebra that can model Euclidean transformations in a struc\-ture-preserving manner}, a distinction sometimes claimed for CGA (\cite{dfm07}, p. 364).  The importance of this result will become more apparent  in the next section,  which turns to a  practical comparison of the two models.  
\section{Comparison: practical issues}
\label{sec:comparison}
Given the same feature set in these algebras, 
we shift our comparison to  more practical considerations for \thetask.  

\subsection{Complexity} 
The point $\vec{x} = (x,y,z) \in \Euc{3}$ receives coordinates $(x,y,z,1)$ in PGA and $(1,x,y,z,\frac12 \| \vec{x} \|^2)$ in CGA.\footnote{This parametrization produces a paraboloid of revolution as null quadric.  To obtain the unit sphere as null quadric, one can rotate by $\frac{\pi}{2}$ in the plane of the two \quot{extra} dimensions to obtain the coordinates $(\frac12 (\| \vec{x} \|^2+1), x, y, z, \frac12 (\| \vec{x} \|^2-1))$.}  The last coordinate is clearly non-linear function of the original ones.  This standard representation is sometimes called the \emph{zero radius sphere} (ZRS) representation of points in CGA.  
Also note, that as a result of having 1 more dimension, CGA also has twice the number of dimensions as PGA.

A more serious effect of the non-linear embedding of this representation, is that  flat euclidean geometric configurations have to be represented and calculated as intersections of linear configurations with the null sphere of $\pclal{n+1}{1}{0}$.  For example, if you want to subdivide a polygon in PGA, linear interpolation will preserve the flatness of the polygon; in the ZRS representation of CGA,  you have to follow linear interpolation with a  projection back onto the null sphere or devise other interpolation methods.  As mentioned above in \Sec{sec:cga}, one can alternatively use the flat representation in CGA, which is essentially the same as PGA.  But to access the distinctive features of CGA (such as the distance function via the inner product) you have then to convert from the flat back to the ZRS representation,  leading to the conclusion that one cannot in this way avoid the consequences of the non-linear embedding. 


In general, \emph{any} computation applied to a geometric primitive in the standard CGA representation risks moving off the null sphere, so potentially each step has to be checked against an error tolerance and corrected.  PGA does not suffer from this difficulty, since all coordinates represent valid euclidean or ideal elements.  Related problems arise with differential equations, as the next section shows.  

\subsection{Numerical analysis and  differential equations} \label{sec:naade} Consider the example of the euclidean equations of motion for a rigid body.   (In PGA, for $n=3$, these are  the biquaternion ODE's given above in \Sec{sec:biq}.)  As with any ODE, the difficulty of solving the ODE is directly connected with how the  space of valid solutions sits inside the full space of possible solutions; the smaller the co-dimension of the former in the latter, the easier the solution process is; see for example the discussion above in \Sec{sec:qmet} regarding the advantages of the quaternion representation of rigid body motion over the matrix one.   
This difficulty is acknowledged in a recent article on 3D euclidean rigid body motion in the conformal model (\cite{lasenby2011}, \S 1.3.1): 
\begin{quote}
... The idea here is to work in an overall space that is two dimensions higher than the base space, using the usual conformal Euclidean setup.  The penalty for doing this, i. e., using a Euclidean setup, is that the number of degrees of freedom is not properly matched to the problem in hand, and \emph{we have to introduce additional Lagrange multipliers to cope with this.} [our italics]
\end{quote}
How do the two algebras compare in this regard? In both, the space of valid solutions for the Euler equations of rigid body motion consists of a euclidean bivector (in the Lie algebra $\mathfrak{e}(3)$) plus a euclidean rotor (in the Lie group $\Eucg{3}$); each is 6D, so the space of valid solutions is 12-dimensional.
In the case of CGA, the  full space of possible solutions is 26-dimensional (10D bivector and 16D even subalgebra).  Here the co-dimension is 14, larger than the solution space itself, forcing the use of Lagrange multipliers.  In contrast, the full  solution space in PGA is 14-dimensional (6D bivector and 8D even subalgebra). 
Normalizing the rotor in PGA to have unit norm brings the solution back onto the valid solution space and  provided reliable results for the extensive simulations in \gTh, Ch. 12), although  the use of a Lagrange multiplier for the same purpose in the PGA case should be investigated.  Given the availability of a fast and reliable PGA solution with minimal need for Lagrange multipliers, the question naturally arises, what advantages does the CGA approach to rigid body mechanics offer in compensation?
\subsection{Kinematics, rigid body mechanics, and classical screw theory}  \label{sec:krbdasct} As we noted above, PGA contains within it the biquaternions and their elegant representation of 3D euclidean kinematics and rigid body mechanics.  This is essentially also the content of the screw theory of Robert Ball \cite{ball00}. As a result, all the features of these theories are included in PGA as a \quot{native} element. 
\cite{hestenes10}, on the other hand, envisions CGA as a means to \quot{rejuvenate} classical screw theory.  It would be worthwhile to compare the two approaches to screw theory with regard to such criteria as compactness of expression, practicality, and comprehensibility.  For example, screw theory is  built upon  line geometry in $\RP{3}$, which in turn is built upon  Pl\"{u}cker coordinates for lines (bivectors).  As noted above, these coordinates are native to PGA, but have to be extracted from the 10-dimensional coordinates of bivectors in CGA using the so-called conformal split (\cite{hestenes10}, \S VII).  

\subsection{Learning curve} \label{sec:learn} The simplicity of the PGA representation leads to considerable savings in explaining the model, a significant advantage when considering the unfamiliarity of the underlying concepts.  Since CGA is embedded in the hyperbolic PGA ($\pdclal{n+1}{1}{0}$), learning the projective model is a natural  step towards understanding the conformal model, but not vice-versa. From a pedagogical point of view, PGA forms a natural transition step between VGA and CGA.  (A  \emph{vector geometric algebra}, or VGA for short, is a geometric algebra whose elements are interpreted as elements of a vector space, rather than projective space. )

\section{Roundness and CGA}
\label{sec:round}
We have established that for \thetask, PGA exhibits a series of practical advantages.  Most of these can be arise from the contrast between, on the one hand, the flat embedding of $\Euc{n}$ in $\RP{n}$ in PGA and, on the other,  the curved embedding of $\Euc{n}$ within  $\RP{n+1}$ in CGA.  We could say that the \quot{roundness} of CGA is a liability when one is restricted to flat primitives.  There is however another \quot{roundness} in CGA that, for some euclidean applications,  compensates for these liabilities: euclidean spheres are represented as points in CGA and can be manipulated on the same level as traditional flat primitives: points, lines, and planes.  The restriction to spheres of radius 0  yields the curved model of $\Euc{n}$ which has formed the basis of the comparison up to now.  Removing the restriction to traditional flat primitives yields a powerful geometric toolkit ideally suited for many euclidean tasks where spheres (or conformal maps, see \Sec{sec:cga}) play an intrinsic role (\cite{dfm07}, Ch. 14).  There are a number of application areas, from optimization to robotics, whose problem settings do exhibit this close connection to sphere (or conformal) geometry.  See, for example, \cite{dorst2014}.

Users choosing between PGA and CGA are therefore advised to carefully weigh the advantages and disadvantages of the \quot{roundness} of CGA in their decision.  On the one side are the advantages of having a direct, powerful representation of spheres; on the other hand, are the disadvantages (discussed above) arising from the embedding of euclidean space \emph{itself} as a sphere in a higher-dimensional projective space (the null sphere of $\pclal{n+1}{1}{0}$). 
Hence, if there is no intrinsic need for sphere geometry, as is the case for applications based on flat geometry,  classical kinematics and rigid body mechanics,  then the disadvantages listed in \Sec{sec:comparison} can be expected to outweigh the advantages; if spheres  form an essential geometric primitive, CGA is probably the right tool for the job.  In between there's much room for further research and development on how the strengths and weaknesses of the two approaches can be optimally combined.

\section{Conclusion}
\label{sec:conclusion}
We have reviewed important concepts from $19^{th}$ century mathematics, and clarified a set of fundamental terms, including \emph{homogeneous model}, \emph{euclidean},  $\R{n}$,  and \emph{duality}, which are key to a correct understanding of how geometric algebra can be applied to doing euclidean geometry. On this basis, we have established that the dual projective geometric algebra $\pdclal{n}{0}{1}$ deserves the title of \quot{standard} homogeneous model for euclidean geometry. 
We have shown that it exhibits all the attractive features with respect to doing euclidean geometry which modern geometric algebra users expect, and is the \emph{smallest} such algebra.  Furthermore, in regard to practical considerations for \thetask, we have found that it exhibits advantages over the higher-dimensional CGA, which  remains the tool of choice for applications making essential use of spheres or conformal maps.  For $n=3$, the most popular case, the fact that $\pdclal{3}{0}{1}$ is built atop the biquaternions, William Clifford's \quot{other} great discovery
 (besides geometric algebra), means that users already familiar with the biquaternions are well-positioned to acquire PGA skills quickly.  

\begin{figure}[h]
  \centering
   \setlength\fboxsep{0pt}\fbox{
    \includegraphics[width=0.6\textwidth]{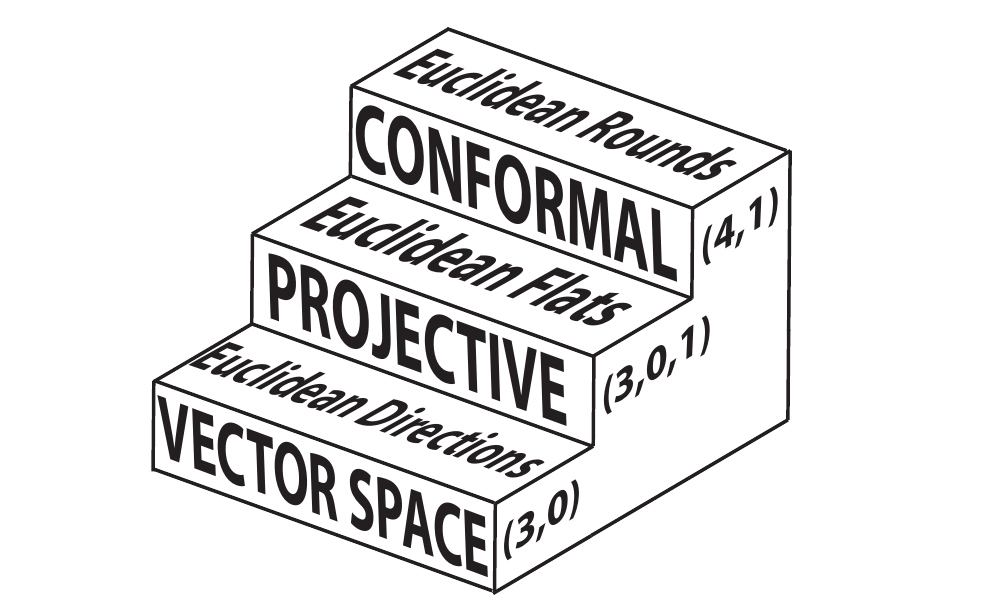}
    }
  \caption{Progression of geometric algebras for euclidean geometry.}
\label{fig:stairs}
\end{figure}

One area of concern for all practitioners of geometric algebra, regardless of specialization, is the slow rate at which geometric algebra has been adopted into the university curriculum (see \cite{ggap2011}, p. vi).  We believe that  the foregoing comparison can make a significant contribution to a solution of this challenge.  PGA and CGA are not directly competitive, any more than automobiles and airplanes are directly competitive.  As the previous section hopefully indicates, each has their proper, interdependent place in the geometric algebra ecosystem. 
 Pedagogically, $\pdclal{n}{0}{1}$ provides a natural stepping-stone (the \quot{automobile}) between the $n$-dimensional vector space algebra $\clal{n}{0}{0}$  (the \quot{bicycle}  of the GA world, aka VGA) and the $(n+2)$-dimensional CGA (the \quot{airplane}).  See \Fig{fig:stairs}. 
 In light of \Sec{sec:learn}, we can expect that PGA will be accessible to a significantly larger pool of  students than currently is the case with CGA.  This will also simplify the teaching of CGA since, as remarked above, CGA  is built on top of PGA.
We suggest that the PGA treatment of euclidean geometry, kinematics, and mechanics is exactly what is needed to solve the GA adoption problem by providing a non-trivial link between VGA and CGA, just as the automobile fits between the bicycle and the airplane. 

Thus, these two algebras exist not in a competitive but a complementary relationship. The nature of this complementarity was already expressed by Johannes Kepler -- perhaps the first scientist to apply mathematics in a modern way to the study of the outer world -- when he wrote (\cite{keplerHM}):  \begin{quote}...God, in his measureless wisdom, selected at the very beginning the straight and the curved, in order with them to imprint into the world
the divinity of the Creator... In this way the most Wise devised the extensive
world, whose whole being is encompassed between the two contrary
principles, the straight and the curved.\end{quote}   Rounding off our article in Kepler's style, we might say that God has given us PGA to understand and to master  the world of the straight, while he has given us CGA to do the same for the world of the round.

\section{Acknowledgements} The author would like to thank the third reviewer for his detailed comments, which led to numerous improvements in the article.

\begin{appendices}
\section{Versors for the euclidean plane.}
  Begin with two normalized 1-vectors $\vec{a}$ and $\vec{b}$  in $\pdclal{2}{0}{1}$, each representing a line in the euclidean plane and assume that they meet in a euclidean point $\vec{P}$, i. e., $\vec{a} \wedge \vec{b} = \vec{P} \neq 0$.  We show  that $R_{\vec{a}}(\vec{b})$, the reflection of the line $\vec{b}$  in the line $\vec{a}$,  is given by $\vec{a}\vec{b}\vec{a}$.  

Using basic facts of elementary geometry, it's not hard to show that $\vec{x} := R_\vec{a}(\vec{b}) $ can be uniquely characterized as the line which passes through $\vec{P}$, the intersection of $\vec{a}$ and $\vec{b}$, and whose oriented angle to $\vec{a}$ is equal but opposite to the angle $\vec{b}$ makes to $\vec{a}$.  The derivation of the signature $(2,0,1)$ in \Sec{sec:caykln} allows us to translate these conditions into the language of the geometric product in $\pdclal{2}{0}{1}$:  $\vec{x} \cdot \vec{a} = \vec{b} \cdot \vec{a}$ and $\vec{x}\wedge \vec{a} = -\vec{b}\wedge\vec{a}$.  First note that $\vec{a}\vec{b}\vec{a}$ is a 1-vector (a line), since the three arguments are not linearly independent, so their wedge is 0.  Is it the desired line?  Substituting $\vec{a}\vec{b}\vec{a}$ for $\vec{x}$ in the first condition and applying symmetry of ($\cdot$) and the normalization condition $\vec{a}^2=1$ yields 
\begin{align*}
\vec{x}\cdot \vec{a} &= \dfrac12 (\vec{a}\vec{b}\vec{a}^2+\vec{a}^2\vec{b}\vec{a}) \\
&=  \dfrac12 (\vec{a}\vec{b}+\vec{b}\vec{a}) \\
& = \vec{a}\cdot\vec{b} = \vec{b} \cdot \vec{a}
\end{align*}
And the second condition proceeds similarly, using  the anti-symmetry of $\wedge$:
\begin{align*}
\vec{x}\wedge \vec{a} &=  \frac12 (\vec{a}\vec{b}\vec{a}^2-\vec{a}^2\vec{b}\vec{a}) \\
&=  \frac12 (\vec{a}\vec{b}-\vec{b}\vec{a}) \\
& = -\vec{b}\wedge\vec{a}
\end{align*}
Hence, $\vec{x}$ fulfills the conditions and  is, therefore,  the reflection of the line $\vec{b}$ in the line $\vec{a}$.  We leave it as an exercise for the reader to verify that this argument remains true when $\vec{a}$ and $\vec{b}$ are parallel (i. e., $\vec{P}$ is ideal).  A further exercise: $\vec{a}\vec{P}\vec{a}$ is the reflection of a euclidean point $\vec{P}$ in the line $\vec{a}$.  The tireless reader can then extend this result to the full euclidean isometry group by applying the well-known result that all isometries can be factored as a sequence of reflections in euclidean lines. Such a sequence of reflections yields, in the algebra, a  versor $\vec{R}$ consisting of the product of the corresponding 1-vectors.  When the 1-vectors are normalized, then so is the resulting versor, which then belongs to the rotor group, and one obtains the sandwich operation associated to this rotor: $\vec{R} \vec{X} \widetilde{\vec{R}}$. Nothing in this proof essentially depends on the dimension $n=2$, it generalizes directly and establishes the claim that the versor representation of the isometry group  in $\pdclal{n}{0}{1}$ works as advertised.

 \end{appendices}


\end{document}

%% file: GunnMac.tex
\usepackage{xspace}
\usepackage{helvet}
\usepackage{courier}
\usepackage{fancyhdr}
\usepackage{type1cm}         
\usepackage{subeqnarray}
\usepackage{graphicx}        
\usepackage{multicol}        
\usepackage[bottom]{footmisc}
\usepackage{amscd}	
\usepackage{enumerate}
\usepackage{amsmath, amsthm, latexsym}
\usepackage{array}
\usepackage{tabularx}
\usepackage{url}
\usepackage{relsize}
\usepackage{ifthen}
\usepackage{xcolor}
\usepackage[newitem,newenum]{paralist}
 
\def\isBeamer{false}

\newcommand{\mychapterskip}{\ifthenelse{\equal{\isMPK}{true}}
{\vspace{-.8in}}{}}

\newcommand{\fullversioncolor}{\color{black!50!blue}}
\newcommand{\versions}[2]{\ifthenelse{\equal{\isfullversion}{true}}{{\fullversioncolor#1}}{#2}}

\renewcommand{\vec}[1]{\mathbf{#1}}
\newcommand{\quot}[1]{``#1''}

\newcommand{\R}[1]{\ensuremath{\mathbb{R}^{#1}\xspace}}
\newcommand{\Euc}[1]{\ensuremath{\mathbf{E}^{#1}\xspace}}
\newcommand{\Eucg}[1]{\ensuremath{E(#1)\xspace}}

\newcommand{\RD}[1]{\ensuremath{(\mathbb{R}^{#1})^{*}\xspace}}
\newcommand\proj[1]{\ensuremath{\mathbf{P}(#1)}\xspace}
\newcommand\RP[1]{\ensuremath{\mathbb{R}{P^{#1}}\xspace}}

\newcommand{\Sph}[1]{\ensuremath{\mathbf{S}^{#1}\xspace}}
\newcommand{\Ell}[1]{\ensuremath{\mathbf{Ell}^{#1}\xspace}}
\newcommand{\Hyp}[1]{\ensuremath{\mathbf{H}^{#1}\xspace}}

\newcommand{\e}[1]{\vec{e}_{#1}}

\newcommand{\eye}{\vec{I}}

\newcommand{\MN}{metric-neutral\xspace}

\newcolumntype{Y}{X}

\newcommand{\bivo}[1]{{\vec{#1}}}

\newcommand{\momo}{\bivo{\Pi}}

\newcommand{\sigo}{\bivo{\Sigma}}

\newcommand{\grass}[1]{\bigwedge{\mathbb{R}^{#1}}}

\newcommand{\pgrass}[1]{\proj{\bigwedge{\mathbb{R}^{#1}}}}
\newcommand{\pdgrass}[1]{\proj{\bigwedge({\mathbb{R}^{#1})^*}}}

\newcommand{\clal}[3]{\mathbb{R}_{#1,#2,#3}\xspace}

\newcommand{\pclal}[3]{\proj{\mathbb{R}_{#1,#2,#3}}\xspace}
\newcommand{\pdclal}[3]{\proj{\mathbb{R}^*_{#1,#2,#3}}\xspace}

\newcommand{\pdclplus}[3]{\proj{\mathbb{R}^{*+}_{#1,#2,#3}}\xspace}
\newcommand{\pclplus}[3]{\proj{\mathbb{R}^{+}_{#1,#2,#3}}\xspace}

\newcommand{\VS}{\ensuremath{\mathsf{V}}}
\newcommand{\VSD}{\ensuremath{\mathsf{V}^{*}}}

\newcommand{\myboldhead}[1]{\vspace{.1in}\hspace{-.35in}\textbf{#1.}}

\newcommand{\app}[1]{\Sec{sec:J}}

\definecolor{mypurple}{rgb}{1,0,1}

\definecolor{mygreen}{rgb}{0, .5, 0}

\newcommand{\Fig}[1]{Fig.~\ref{#1}}

\newcommand{\Sec}[1]{Sect.~\ref{#1}}

\ifthenelse{\equal{\isBeamer}{true}}
{

}
{

}

\newcommand{\mycorrection}[1]{}

